\documentclass[12pt]{article}
\usepackage{amsmath,amsfonts,amssymb}
\usepackage{graphicx,color}

\usepackage{theorem}

\theorembodyfont{\upshape}

\newcommand{\ZZ}{{\mathbb Z}}

\newcommand{\CC}{{\mathbb C}}
\newcommand{\NN}{{\mathbb N}}

\def\g{{\mathfrak{g}}}

\def\hh{{\mathfrak{h}}}

\def\mm{{\mathfrak{m}}}

\def\O{{\cal{O}}}

\numberwithin{equation}{section}
\begin{document}

\begin{center}
{\bf Generalized Casimir Operators}\\[5mm]
{\bf S. Eswara Rao}\\[5mm]
\end{center}

\begin{abstract}
Let $\g$ be symmetrizable Kac-Moody Lie algebra.  In this paper we describe
a new class of central operators generalising the Casimir operator.  We also
prove some properties of these operators and show that these operators move
highest weight vectors to new highest weight vectors.\\

MSC: Primarty 17B65, secondary 17B10, 17B70, 17B69.\\
\end{abstract}

\section*{Introduction} Let $\g$ be symmetrizable Kac-Moody Lie algebra and
$A$ be a commutative associative algebra with unit.  Then
$\g \otimes A$ is naturally Lie algebra and $\g$ is a subalgebra.  We 
consider a certain
category $\O$ of $\g\otimes A$ modules. See Definition 1.3

We now construct a class of operators $\Omega (a,b), a,b \in A $ which act
on modules in $\O$ and commutes with $\g$ action. Such operators are called
central operators.  These operators are variations of Casimir operator and
in fact $\Omega (1,1)$ is the Casimir operator. It is  well known that Casimir 
operator acts as a scalar on $\g$ highest weight vectors.  Where as our 
central operators move one $\g$ highest weight vector to 
another most often.  This way if we know one highest weight vector by 
applying our central operators we can produce more highest weight vectors of
the same weight.\\

The idea of these central operators was born in trying to understand
evaluation modules.  We will explain this in the simplest case 
$A=   {\CC}[t,t^{-1}]$. Let $V(\lambda_i), \ 1\leq i \leq n $, be an 
irreducible integrable highest weight module for $\g$. Then the tensor 
product module $V=\otimes V(\lambda_i)$ is known to be completely 
reducible as $\g$-module. $V$ can be made into $\g\otimes A$ module by 
evaluating at distinct points (see 3.3) and is called evaluation module for
$\g \otimes A$. In this case there are special central operators denoted
by $\Omega (l,k), 1 \leq l, \ k\leq n$, which act only on the $l$th and $k$th
factors. In fact it is a Casimir operator acting on $V(\lambda _l)
\otimes V (\lambda_k)$ and the identity on the rest of the factors. Linear
combinations of $\Omega(l,k), 1 \leq l,  k \leq n$ exhaust all our central
operators in the evaluation module case.  We have defined highest weight modules 
$V(\psi)$ for $\g \otimes A$ and all evaluation modules are highest weight modules.  But there are many more highest weight modules which are
not evaluation modules. We do not know how these central operators act on
highest weight modules.  When $\g$ is a simple finite dimensional Lie algebra, the decomposition $V$ as a  $\g$ module is a classical open problems.
There are several results available for $n=2$.  See [KU1] and references
therein. But it looks like not much is known for $n\geq 3$ and here our 
central operators are very effective. We work out some examples 
(Examples (3.11) and (3.12)) and note that in these examples that the space 
spanned by repeatedly applying our central operators on a single highest weight
vector gives the whole highest weight space of that weight.  This will not be
 true in general. For example in the case $n=2$ all our central operators are 
scalars on $\g$-highest weight vectors and so not very interesting.

In the last section we consider $\g=gl_N$.   In this case we have more central
operators.  It is known that the center of $U(\g)$ is finitely generated as 
an algebra.  In fact for every positive integer $k$ there is the  $T_k$ (called
Gelfand invariant) in the center of $U(\g)$ and $T_1,T_2,\ldots T_N$ generate
the center of $U(\g)$ as an algebra.  Now for each $k$ we define a class of 
Central Operators (depending on $A$).  See (4.2) and Proposition (4.4).\\

We will now write down these Central Operators explicitly in the evalution 
module case.  As earlier these operators are independent of $A$ in the
case of evalution modules.\\

Let $E_{ij}$ be the standard basis of $\g$.  Recall that the Gelfand invariant
$$
T_k=\sum_{(i_l,\ldots i_k)} E_{i_1i_2} E_{i_2 i_3}\ldots
E_{i_k i_1}.
$$
Fix a positive integer $n$ and let $V_1,V_2,\ldots V_n$ be irreducible
finite dimensional module for $\g$. Consider
$$
V =V_1\otimes V_2 \otimes\ldots\otimes V_n
$$
Let $j$ be such that $1\leq j\leq n$ and define the operator 
$E_{i_1i_2}(P_j)$ on $V$ as $E_{i_1 i_2}$ acting only on the $j$th factor
of $V$.  
Now define for $1\leq j_i\leq n$
$$T_k(P_{j_1},\ldots, P_{j_k}) = \sum_{(i_1,\ldots, i_k)} E_{i_1 i_2}(P_{j_1})
\ldots  E_{i_k i_1}(P_{j_k}).$$

We prove all these operators are central that is, they commute with 
the $\g$ action 
on $V$. Further the original Gelfand invariant
$$
T_k =\sum_{l \leq j_1,\ldots,j_k \leq n}
 T_k(P_{j_1},\ldots, P_{j_k}) 
$$

It is well  known that $V$ decompose into irreducible finite dimensional
$\g$-modules.  It is also known that each $T_k$ acts as a scalar on any 
$\g$ isotypic component of $V$. Whereas the operator $T_k(P_{j_1},\ldots, P_{j_k})$ 
does not act as scalars on these $\g$-components.  They take one $\g$ highest
weight vector to a new $\g$ highest weight vector most often.  These operators  will be 
greatly useful for finding highest weight vectors once we know one highest
weight vector.  We will now write down one such operator explicitly.  Take
$k=4$ and $n=4$.  Let $w_i \in V_i$\\
$$
T_4(P_2,P_1,P_2,P_3)(w_1 \otimes w_2 \otimes w_3 \otimes w_4) =
\sum_{(i_1,i_2,i_3,i_4)} E_{i_2 i_3} w_1 \otimes E_{i_1 i_2} E_{i_3 i_4} w_2
\otimes E_{i_4 i_1} w_3 \otimes w_4
$$
Let  $T$ be non-commutative associative algebra generated by $T_k(P_{j_1},
\ldots, P_{j_k})$ for all $1\leq j_1,\ldots, j_k\leq n$ and for all $k> 0$.
Then in Theorem 4.8, we note that each isotypic component is irreducible
for $U(\g)\otimes T$.
Recall that $T$ is an algebra generated by Central Operators $T_k(a_1,a_2,
\ldots, a_k)$.  In particular it contains finite products of such operators.  
In the last section we will give a spanning set and will avoid products.  We
will define certain twisted operators which are again central and prove that
it is spanning set for $T$.
\paragraph*{Section 1}

Throughout the paper  all vector spaces and tensor products are over complex
numbers ${\CC}$. $U$ always denotes the universal enveloping algebra of a 
Lie-algebra

\paragraph*{(1.1)} Let $\g$ be a symmetrizable Kac-Moody Lie algebra.  Let  (,) be a 
non-degenerate invariant symmetric bilinear  form on $\g$.  Let $\hh$  be a  
Cartan subalgebra. Let $\{ \alpha_1,\ldots \alpha_l\}$ and $\{ \alpha^\vee_1,
\ldots, \alpha^\vee_l\}$  be roots and coroots of $\g$. Let $\Delta$ and 
$\Delta^+$ be roots and positive roots of $\g$.

Let 
$$ \g=\bigoplus_{\alpha \in \Delta} \g_\alpha \oplus \hh
$$
be the root space decomposition of $\g$.  See Kac book $[K]$ for more details.

\paragraph*{(1.2)} Let $A$ be a  commutative associative algebra with unit.
Denote $\g(A)=\g\otimes A$ with obvious Lie bracket.  For any vector space
$V$ denote $V(A)=V\otimes A$.  Let $\g=N^+ \oplus \hh\oplus N^{-}$ be the
 standard
triangular decomposition.  Then $\g(A)=N^+(A) \oplus \hh(A) \oplus 
N^{-} (A)$ is a triangular decomposition for $\g(A)$.  For 
$\alpha \in \Delta^+$ 
define  $ht\ \alpha=\sum n_i$ where
$\alpha =\sum n_i \alpha_i$.  Note that $\g \simeq \g \otimes 1$.

\paragraph*{(1.3)} Definition: $A$ module $V$ of $\g (A)$ is said to be in the  category 
$\O$ if the following holds 

(a) $V$ is a weight module for $\g(A)$ with respect to the Cartan subalgebra 
$\hh$ and has finite dimensional weight spaces.
(b) For every  $v$ in $V$ and $a \in A$ we have $(X_\alpha \otimes a) v = 
0$ for $ht\alpha\gg 0$ and $\alpha \in \Delta ^+$ and $X_\alpha \in 
\g_\alpha$.

\paragraph*{(1.4)} We will now produce a class  of irreducible $\g(A)$ modules which are 
in $\O$. Let $\psi:\hh(A) \rightarrow {\CC}$ be any linear map. Consider
the one
 dimensional vector  space ${\CC} v$ which is $N^+(A) \oplus \hh(A)$ module where $\hh(A)$ acts  via $\psi$ and $N^+(A)$ acts trivially.  Now consider the Verma  module.
$$
M(\psi) =U(\g(A)) \bigotimes_{N^+(A) +\hh(A)}{\CC v}.$$
By standard arguments  we see that $M(\psi)$ has an unique irreducible quotient denoted by $V(\psi)$. Note that when $A$ is infinite dimensional $M(\psi)$ does not have finite dimensional  weight spaces.  $V(\psi)$ may have finite dimensional weight spaces depending on $\psi$.

Let $\g'=[\g,\g]$ and let  $\hh' =\g\cap \hh$.  Let $\hh^{''}$ be
 any vector space such that $\hh=\hh' \oplus \hh^{''}$.  See $[K]$ for more 
details.  Let $\stackrel{\sim}{\g} =\g' (A) \oplus \hh^{''}$. Lie algebra
 $\stackrel{\sim}{\g}$
was originally considered in [E3] and module theory is developed for the 
special case  where $A$ is a Laurent polynomial algebra in several commutating 
variables. They have been generalised for any $A$ in [EB].

\paragraph*{(1.5) Lemma:} $V(\psi)$ is irreducible as 
$\stackrel{\sim}{\g}$ module.

\paragraph*{Proof} First note that $U(\stackrel{\sim}{\g}) v=V(\psi)$ as the 
additional
space $\hh^{''} \otimes A$ acts a scalars on $v$.  Suppose $W$ is a 
$\stackrel{\sim}{\g}$ submodule of $V(\psi)$. Let $ w \in W$ be a weight  vector of 
maximal height.  Then clearly $w$ is a highest weight vector in the sense that 
$(\g_\alpha \otimes A) w=0$ for all $\alpha \in \Delta^+$.  But 
$V(\psi)$ does not have highest weight vectors except the multiples of $v$.  Thus $w=v$ upto scalar.  This proves $W=V(\psi)$.  Lemma is proved.

Since $V(\psi)$ is an irreducible $\stackrel{\sim}{\g}$ - module,  we can 
use results from [EB].

\paragraph*{(1.6) Proposition (Prop. 2.4 and Lemma 2.3, [EB])}

$V(\psi)$ has finite dimensional weight spaces if and only if there exists a 
co-finite ideal $I$ of  $A$ such that $\g'\otimes I \cdot V(\psi)=0$.

\paragraph*{(1.7)} Such $V(\psi) \in {\cal O}$.
\paragraph*{(1.8)} There exists a special class of co-finite ideals. Fix a 
positive 
integer $n$. Let $\mm_i \ 1 \leq i \leq n$, be distinct maximal ideals of
$A$ and assume $A$  is finitely generated.
Because of the assumptions on $A$  we know that $A/ \mm_i \cong 
{\CC}$.
Consider the co-finite ideal $I=\cap \mm_i$.  Then by Chinese Reminder  
Theorem we have $A/I \cong \oplus{\CC}$ so that 
$\g\otimes A/I \cong \oplus (\g\otimes A/\mm_i) \cong \oplus \g$.  
For each $i$ let $V(\lambda_i)$ be an irreducible highest weight  module for 
$\g$ with highest weight vector $v_i$ and highest weight $\lambda_i$. Then 
$V=\bigotimes^n_{i=1} V(\lambda_i)$ is a  
irreducible  $\g(A)$ module via the surjective map $\Pi:\g(A)\rightarrow  
\oplus \g$. Note that the space $\oplus \hh$ acts as scalars on $v=v_1 
\otimes \ldots \otimes v_n$ and now consider the  surjective map $\hh(A)
\rightarrow \oplus \hh$. Let the corresponding map from $\hh(A) \rightarrow 
{\CC}$ by $\psi$.  Then it is easy to see that $V(\psi)\cong V$ as 
$\g(A)$ -modules.

\paragraph*{(1.9)} Such modules $V(\psi)$ are called evaluation modules.

Several generalisation of evaluation modules are considered in the literature. 
See [NS] and references there in.  For super case see [S].

\paragraph*{2. Section: Central Operators }
\paragraph*{(2.1)} We will first recall a certain classical problem in Lie 
theory. We assume $\g$ is simple finite dimensional Lie algebra.  Let 
$V_1,V_2,\ldots V_n$ be irreducible finite dimensional $\g$-modules.  Then 
$V=\otimes  V_i$ be the tensor product module for $\g$.  It is well known 
that $V$ is completely reducible as $\g$-module.

\paragraph*{(2.2) Open problem : Which $\g$ modules occur in $V$ and 
with what multiplicity?}

\vspace*{0,1in}
There are several results available in the literature and most often for
$n=2$.  See [KS1] and references there in. We will now define a class of 
operators, which generalise Casimir operator, acts on the tensor product 
module $V$ and commutes with $\g$.  The main property of our operators, when
applied on a $\g$ highest weight vector, produces a new highest weight 
vector.  Whereas the Casimir operator acts as scalar.  We will define our 
operator in the generality of symmetrizable Kac-Moody Lie algebra and they 
are central operators in the following sense.

\paragraph*{(2.3) Central Operators:} A linear operator acting on objects
 of $\O$ is called central operator if it commutes with $\g$ action.

We will now closely follow Chapter 2 of Kac book [K].  Let  $\hh^*$ be the 
dual of the Cartan subalgebra $\hh$ and denote the non-degenerate symmetric 
bilinear form as (,).  We have an isomorphisms

$\nu: \hh\rightarrow \hh^*$ defined by 
$$
< \nu (h),h_1 >= \nu (h) (h_1)= (h, h_1)\\
$$

Let $ \rho\in \hh^*$ be such that $(\rho, \alpha_i) =\frac{1}{2}(\alpha_i,
\alpha_i), 1\leq i \leq l$. 
Let $\{ e^j_\alpha\}$ be a basis of $\g_\alpha$ and let $\{ e^j_{-\alpha}\}$ 
be the dual basis.  Let $u_1, u_2,\ldots, u_l$ be a basis of $\hh$ and let 
$u^1,u^2,\ldots, u^l$ be the dual basis.  Let $x(a) =x\otimes a, \ x\in \g,
 \ a\in A$.  For $a,b \in A$, Define
$$
\Omega_{a,b}= \sum_{\alpha \in \Delta^+} \sum_j e^j_{-\alpha} (a) 
e^j_\alpha(b)
$$
Now define the operator
$$\Omega(a,b) =2\nu^{-1} (\rho) (ab) + \sum_i u^i(a) u_i(b)+\Omega_{a,b} + \Omega_{b,a}
\leqno{\bf (2.4)}$$
Certainly $\Omega(a.b)$ is infinite sum and sits inside some completion of 
$U(\g(A))$. But $\Omega(a,b)$ is  locally finite on any $V$ in $\O$. 
(Note that it preserve the weight spaces).
In the sense, given a $v$ in $V\in \O$ then $\Omega (a.b) v$ is a finite sum.
We also note that $\Omega(a,b)$ is linear in both variable.  That is 
$\Omega(\lambda_1 a_1 +\lambda_2 a_2, b) =\lambda_1\Omega (a_1,b) + 
\lambda_2\Omega(a_2, b)$ and the same thing is true in $b$ also,  for $a_1, a_2, b\in 
A,\  \lambda_1, \lambda_2 \in {\CC}$.

\paragraph*{(2.5) Theorem:} For  $a,b \in A, \ \Omega (a,b)$ is a central
operator on modules in $O$.\\
We first prove some Lemmas.

\paragraph*{(2.6) Lemma:} Let $ a,b \in A$ and let $\alpha, \beta \in 
\Delta$. Let $z \in \g_{\beta-\alpha}$. Then  we have

$$
\sum_s e^s_{-\alpha} (a) [z, e^s_{\alpha} (b)]
= \sum_s [e^s_{-\beta} (a), z] e^s_\beta (b) \in U (\g(A))
$$
\paragraph*{Proof} First recall the following  Lemma (2.4) from $[K]$.

$$
\sum_s e^s_{-\alpha} (1) \otimes [z, e^s_\alpha (1)] 
=\sum_s[e^s_{-\beta}(1), z] \otimes e^s_\beta (1) \in \g \otimes \g \leqno{(2.7)}
$$
Now consider the following $\g$-module homomorphism from \\
$$
\g \otimes \g \, \,\, \text{to}\,\,\,  U(\g(A))
$$
sending $X\otimes Y$ to $ X(a) Y(b)$.  Applying the $\g$-modules 
homomorphism to 2.7, Lemma 2.6 follows.

\paragraph*{(2.8) Lemma } Let $a,b \in A$

$$
\begin{array}{lll}

(1) & [\Omega_{a,b}, e_{\alpha_i}]=-\nu^{-1} (\alpha_i) (a)
 e_{\alpha_i}(b)\\
(2) & [\Omega_{b,a}, e_{\alpha_i}]=-\nu^{-1} (\alpha_i) 
(b) e_{\alpha_i}(a)
\end{array}
$$

\paragraph*{Proof} Proof is similar to the proof of theorem 2.6(a) 
of [K]. 
See the second part on pager 22.  We need to use Corollary (2.10).

\paragraph*{(2.9) Lemma}
$$
[\sum u^j(a) u_j(b), e_{\alpha_i}]=\nu^{-1} (\alpha_i)(a) e_{\alpha_i}(b)
+ e_{\alpha_i} (a)  \nu^{-1} (\alpha_i)(b)
$$
Direct checking using 2.5.3 of [K]. Also use the fact that
$\nu$ preserves the bilinear form on $\hh$ and $\hh^*$.

\paragraph*{(2.10) Lemma} For $\alpha \in \Delta$
$$
\begin{array}{ll}
(a) & \alpha (\nu^{-1} (\rho)) =(\rho, \alpha)\\
(b) & \alpha(\nu^{-1} (\alpha)) =(\alpha, \alpha)
\end{array}
$$
Just use the 2.5.3  of [K] 

\section*{Proof of Theorem (2.5)}
From the above Lemma we see that
$$
\begin{array}{lll}
&[\Omega(a,b), e_{\alpha_i}] =[2\nu^{-1}(\rho)(ab), e_{\alpha_i}]\\
&+\nu^{-1}(\alpha_i)(a)e_{\alpha_i}(b) +e_{\alpha_i}(a)\nu^{-1}
(\alpha_i)(b)\\
&-\nu^{-1}(\alpha_i)(a) e_{\alpha_i}(b) -\nu^{-1}(\alpha_i)(b)
e_{\alpha_i}(a)
\end{array}
$$
Note the first term is equal to
$$
2\alpha_i(\nu^{-1}(\rho))e_{\alpha_i}(ab)=2(\rho, \alpha_i) e_{\alpha_i}(ab)
$$
Also note that
$$
\begin{array}{l}
e_{\alpha_i}(a)\nu^{-1} (\alpha_i)(b)-\nu^{-1} (\alpha_i)(b) e_{\alpha_i}(a)
= -\alpha_i (\nu^{-1}(\alpha_i))e_{\alpha_i}(ab)\\
=-(\alpha_i, \alpha_i)e_{\alpha_i}(ab)
=-2(\rho, \alpha_i) e_{\alpha_i}(ab)
\end{array}
$$
Now it is easy to see
$$
[\Omega (a,b), e_{\alpha_i}]=0
$$
In a similar way we see that $[\Omega(a,b), e_{-\alpha_i}]=0$.  Since
$\Omega(a,b)$ zero weight operator it commutes with $\hh$. As $\Omega (a,b)$ 
commutes with all generators of $\g$, it commutes with $\g$.  This completes 
the proof of the Theorem.

\paragraph*{(2.11) Remark:} Theorem (2.5) holds in the generality of 
Borcherds-Kac-Moody super algebras (BKM).  Note that BKM super algebra admits a
unique (upto scalar) non-degenerate, super invariant and super symmetric 
billinear form.  See Theorem 18.4.2 of [MU] for the special case 
$A= \CC.$ See [SN] and [W] for definitions of BKM super algebras.

\paragraph{3. Section}
\paragraph*{(3.1)} Throughout this section we assume $A={\CC}[t, t^{-1}$] a
Laurent polynomial algebra. For any vector space $V$ we denote 
$L(V)=V\otimes A$.

In this section we give three examples to indicate the importance of our 
operators.  We work with evaluation modules and they have been mentioned
in (1.9).  In our case they can be made more explicit. First we will simplify
our central operators on evaluation modules.

We first recall evaluation modules in the context of ${\CC} [t, t^{-1}]$.
See [E1], [E2] and [E3] for some classification results.

\paragraph*{(3.2)} Let $\g$ be a symmetrizable Kac-Moody Lie algebra and $\hh$
be a Cartan subalgebra.  Fix a positive integer $n$ and let $a_1, a_2, \ldots
a_n$ be  non-zero distinct complex numbers.

Let $V(\lambda_1), V(\lambda_2)\ldots V(\lambda_n)$ be irreducible highest
modules for $\g$. Let $v_1, v_2,\ldots. v_n$ be the corresponding highest 
weight vectors.\\
Let $\underline{\lambda} =(\lambda_1, \lambda_2 \cdots \lambda_n), \ 
\underline{a}=(a_1, \ldots, a_n)$\\
Let $V(\underline{\lambda},\underline{a})=\otimes^n_{i=1} V(\lambda_i)$.\\
Define a $L(\g)$ module structure on $V(\underline{\lambda}, \underline{a})$

\paragraph*{(3.3)} $X \otimes t^k (w_1 \otimes\ldots \otimes w_n) =
\sum a^k_i w_1
\otimes \ldots X w_i \otimes \ldots \otimes w_n$ for $X \in \g, \ k \in\ZZ$
and $w_i \in V (\lambda_i)$.  It can easily checked to be $L(\g)$-module. 
We will now indicate another way of seeing this.
Consider the Lie-algebra map

\paragraph*{(3.4)}
$$
\begin{array}{lll}
\Pi(\underline{a}): L(\g)\rightarrow \oplus \g (n \ copies)\\
\Pi (\underline{a}) (X \otimes t^k) =(a^k_1   X, \ldots a^k_n X)
\end{array}
$$
It is standard fact that $\Pi(\underline{a})$ is surjective. See [E3]
\paragraph*{(3.5) Claim:} $V(\underline{\lambda},\underline{a})$ is an 
irreducible as $L(\g)$-module. First note that $V$ is an irreducible module 
for $\oplus\g (n$ copies). Now using the surjective map $\Pi\underline{(a)},
 V(\underline{\lambda},\underline{a})$ becomes $L(\g)$-module and one can
check that this is precisely one given at (3.3). This proves the claim.

Consider $\psi(h\otimes t^k)=\sum a_i^k\lambda_i(h)$ which is linear map from 
$L(\hh)$ to ${\CC}$. Recall we have defined an irreducible  module 
$V(\psi)$ in (1.4).  It is easy to see that $V(\psi) \cong 
V(\underline{\lambda},\underline{a})$ as $L(\g)$-modules.  We will give 
another proof that $\Omega(a,b)$ are central operators.

\paragraph*{(3.6)} Let
$$
\begin{array}{llll}
P(t)  &=& \Pi^n_{i=1} (t-a_i)\\
P_i(t) &=& \frac{\Pi_{i\neq j} (t-a_j)}{\Pi_{i\neq j} (a_i-a_j)}
\end{array}
$$
It is easy to see
\paragraph*{(3.6.1)} $P_i(a_j)=\delta_{ij}$
\paragraph*{(3.6.2)} $\sum P_i(t) =1$

We note that $\g\otimes P(t) (V(\underline{\lambda},\underline{a})=0$ as it
is an evaluation module and $P(a_i)=0$ for all i. Let $I$ be an ideal 
generated
by $P(t)$ and $I$ is a co-finite in $A$. Further $\g\otimes I$
$V(\underline{\lambda},\underline{a})=0$.  Further we note that ker 
$\Pi =\g\otimes I$ (See 3.4). Now it is clear that 
$\Omega (a,b)$ is zero on $V(\underline{\lambda},\underline{a})$ if either 
$a\in I $ or $b\in I$.  We also have $P_i(t) {\not\in}
I$ and is easy to check that the image of $P_i(t),\ 1 \leq i \leq n$ in $A/I$
form a basis for $A/I$.

Thus to consider $\Omega(a,b)$,  we can assume $a$ and $b$ are linear 
combinations of $P_i(t)$.

\paragraph*{(3.7)}These polynomials $P_i(t)$ are very special. For example 
$$
X \otimes P_i(t) (w_1\otimes\ldots \otimes w_n) =w_1 \otimes \ldots 
X w_i\otimes \ldots w_n
$$
Where $X \in\g, \ w_i\in V(\lambda_i)$. So $X\otimes P_i(t)$ acts only on the 
factor $V(\lambda_i)$. This means $\Omega(P_i(t), P_i(t))$ acts only on 
$i$the factor and it can be seen to be the classical Casimir operator 
acting on $V(\lambda_i)$.  In particular it is a central operator. Similarly 
$X\otimes (P_i(t) +P_j(t))(w_1\otimes \ldots w_n)= w_1\otimes \ldots 
Xw_i\otimes \ldots w_n + w_1 \otimes \ldots Xw_j\otimes \ldots \otimes w_n$.

So it will act on $i$th and $j$th factor.  Then the operator 
$\Omega (P_i(t) +P_j(t), \ P_i(t) +P_j(t))$ acts only on $i$th and 
$j$th factor.

It can be readily seen to be classical Casimir operator acting on 
$V(\lambda_i)\otimes V(\lambda_j)$.  It is a central operator. Now we have 
different proof that $\Omega(P_i(t), P_j(t))$ is a central operator. Now from 
above we know that $\Omega(a,b)$ is linear combination 
of $\Omega (P_i(t), P_j(t))$.  Thus it is another proof that $\Omega(a,b)$ are
central operators on an evaluation module.

\paragraph*{(3.8)} We will now digress a little to explain evaluation modules
 in the
 context of finitely generated commutative associated algebra $A$ with unit
 1.
See (1.8) and (1.9) where we have considered evaluation modules of 
$\g\otimes A$,  Let $\mm_1,\ldots \mm_n$ be distinct maximal ideals and we have
 $A/m_i\cong {\CC}$.  We also have surjective map
 $\Pi: A\rightarrow \oplus A/m_i =\oplus {\CC} (n$ copies).  Consider
 $z_i=(0,\ldots 1,\ldots 0)\in \oplus {\CC}$. Let $P_i\in A $ such that 
$\Pi(P_i) =z_i$.  Then clearly $P_1,\ldots, P_n$ is a basis of $A$ mod $I$
 where 
$I=\cap \mm_i$.  Consider the corresponding evaluation module $V(\psi)$ as define
in (1.8). It  is clear that $\g\otimes I \cdot V(\psi)=0$. As explained in 
(3.7) the operators $\Omega(a,b)$ are linear combination of 
$\Omega (P_i, P_j)$. Again $\Omega (P_i, P_i)$ and $\Omega (P_i + P_j, 
\ P_i + P_j)$  are standard  Casimir operator acting  on $V(\lambda_i)$ and 
$V(\lambda_i)\otimes V(\lambda_j)$.  Certainly each of them are central and 
hence $\Omega(a,b)$ is a central operator. It is another proof $\Omega(a,b)$ is central.

\paragraph*{(3.9)} We note that, for an evaluation module, we do not get any new 
central operators for general $A$.  It is sufficient to take 
${\CC}[t, t^{-1}]$.

\paragraph*{(3.10) Remark:} Even though the operators on evaluation module case,
looks familiar we do not have any evidence that they have been considered by 
other authors.  These operators applied to highest weight vector produce new 
highest weight vectors most often. We will explain this with some examples.

\paragraph*{(3.11) Example:} Let $\g$ be any symmetrizabe Kac-Moody Lie
algebra with  the standard form (,). Fix a  positive integer $n$ and 
consider 
$V(\lambda_1),\ldots, V(\lambda_n)$ irreducible highest weight modules for 
$\g$ with highest weight vectors $v_1,\ldots, v_n$  and highest weights 
$\lambda_1, \lambda_2, \ldots \lambda_n$ which we assume to be dominant 
integral.  We know that $V=\otimes^n_{i=1} V(\lambda_i)$ is completely 
reducible $\g$-module. Put $\lambda=\sum\lambda_i$.  Let $V=
\oplus_{\beta\geq 0} V_{\lambda-\beta}$ be weight space decomposition. 
Denote $V^+_\mu$ be the $\g$-highest weight vectors of weight $\mu$.

Let $\alpha_1^\vee,\ldots, \alpha^\vee_l$ be the co-roots. Fix $j$ and 
assume 
$\lambda_i(\alpha_j^\vee)=m_i\geq 1$.  This means 
$e_{-\alpha_j} v_i \neq 0 \forall i$.  Let $w_k =v_1\otimes 
\cdots e_{-\alpha_j} v_k \otimes \ldots \otimes v_n$. Let $z_{k,l} =
m_l w_k -m_k w_l \in 
V_{\lambda-\alpha_j}$.  It is direct checking that $z_{k,l} \in 
V ^+_{\lambda-\alpha_j}$.  We can see that dim $V_{\lambda-\alpha_j} =n$ and 
dim $V^+_{\lambda-\alpha_j}=n-1$.  It is easy to see that $z_{1,2}\ldots,
 z_{1,3}\ldots, z_{l,n}$ are linearly independent and $n-1$ in number.  
Thus it
is a basis for $V^+_{\lambda-\alpha_j}$ .  Since we are working with dual 
basis in the definition of central operators we see that
$$
[e_{\alpha_j}, e_{-\alpha_j}]=\frac{(\alpha_j,\alpha_j)}{2} 
 \alpha^\vee_j$$
Recall the operator $\Omega(a,b)$ and for simplicity let $\Omega (l,k)
 =\Omega (P_l(t), P_k(t))$ for  fixed $l\neq k$.

The following are direct calculation.

\paragraph*{(3.11.1)} (a) $i \not\in \{l,k\}, \ \Omega (l,k) w_i=(\lambda_l,
 \lambda_k) w_i$.\\
(b) $\Omega (l,k) w_k=(\frac{\alpha_j, \alpha_j}{2})m_k w_l + (\lambda_k -\alpha_j, 
\lambda_l)w_k$.

We will now calculate the action of the operators on highest weight vector.  
The following is again direct calculation using 3.11.1.

\paragraph*{(3.11.2)}(a)  Let $p\neq q, \ p,q\not\in \{l.k\}, \Omega 
(l,k) z_{p,q} =(\lambda_l,\lambda_k) z_{p,q}$.\\
(b) $\Omega (l,k) z_{l,k} = ((\frac{(\alpha_j,\alpha_j)}{2} (m_l+ m_k))
-(\lambda_l,\lambda_k)) z_{k,l}$.\\
(c) $q\neq l, \ k\neq q, \ \Omega (l,k) z_{k,q} =(\lambda_l,\lambda_k)z_{k,q} 
-  m_q \frac{(\alpha_j,\alpha_j)}{2} z_{k,l}$.

\paragraph*{3.11.3 Remark} Fix $k\neq l$.  Then applying central operators 
repeatedly on $z_{k,l}$ we get the whole space $V^+_{\lambda-\alpha_j}$

\paragraph*{(3.12) Example:}
Let $\g$ be a symmetric Kac-Moody Lie-algebra. Let $\lambda_1, \lambda_2,
\ldots \lambda_n$ be dominant integral weights. Let $V(\lambda_1),\ldots,
V(\lambda_n)$ be irreducible integrable highest weight modules  with highest
weight vectors $v_1, v_2,\ldots, v_n$.  Let $V=\bigotimes^n_{i=1} 
V(\lambda_i)$ and  let $\lambda =\sum \lambda_i$. Let $v = v_1 \otimes \ldots
\otimes v_n$. Let $\{ \alpha^\vee_1, \alpha^\vee_2, \ldots\alpha^\vee_l
\}$ be co-roots.
Since we are assuming $\g$ to be symmetric we have

\paragraph*{3.12.1} $(\lambda_i,\alpha_j)=\lambda_i(\alpha_j^\vee)$ and $(\alpha_j, 
\alpha_j)=2$. We fix $j$.
\paragraph*{3.12.2} We also assume $m_i =(\lambda_i, \alpha_j) \geq 2$ for all
$i$. This means $e^2_{-\alpha_j} v_i\neq 0$.

Let $V=\bigoplus_{\beta\geq 0} V_{\lambda-\beta}$ be the weight space 
decomposition. $V_\mu^+$ be the space of  $\g$-highest weight vectors. Let
$k\neq l$.
\paragraph*{3.12.3} Let 
$$
\begin{array}{lll}
z_{k,l} &=&v_1\otimes \ldots e_{-\alpha_j} v_k\otimes
\ldots \otimes e_{-\alpha_j} v_l\otimes\ldots v_n.\\
z_k&=& v_1\otimes \ldots e^2_{-\alpha_j} v_k \otimes \ldots \otimes v_n
\end{array}
$$
So that $z_{k,l}, z_k \in V_{\lambda-2\alpha_j}$

\paragraph*{3.12.4} Let 
$$
A_{k,l}=2(m_k-1)(m_l-1) z_{k,l}-(m_k-1) m_k z_l-(m_l-1) m_l z_k
$$

It is direct checking that $A_{k,l} \in V^+_{\lambda-2\alpha_j}$.  Note that 
$A_{k,l}=A_{l,k}$.

The following is easy to see 

\paragraph*{3.12.5}

$$
\begin{array}{lll}
(a)\dim \ V_{\lambda-2 \alpha_j} = \begin{pmatrix}{n}\\{2}\end{pmatrix} +n \\
(b)\dim V^+_{\lambda-2 \alpha_j} =\begin{pmatrix}{n}\\{2}\end{pmatrix} \\
(c) \# \{A_{k,l}, k\neq l\} =\begin{pmatrix}{n}\\{2}\end{pmatrix}
\end{array}
$$
and they form a basis for $V^+_{\lambda-2 \alpha_j}$

The following which gives a formula how our operators act on 
$V^+_{\lambda-2 \alpha_j}$. As earlier let $\Omega(k.l):= \Omega(P_k, P_l)$.
Let $k\neq l$. 

\paragraph*{3.12.6} (a) $  p\neq q,
p,q \not\in \{k,l\}, \Omega (p,q) A_{k,l} =
(\lambda_p, \lambda_q) A_{k,l}$\\
(b) $q=k, \ p\neq l, \ \Omega (p,q) A_{ql} =(\lambda_p, \lambda_q -\alpha_j)
A_{ql} - \frac{(m_l-1)m_l}{(m_p-1)} A_{p,q} + \frac{(m_q-1)m_q}{(m_p-1)}
A_{p,l}$\\
(c) $q=k, \ p=l,\ \Omega(p,q) A_{p,q} =(\lambda_p -\alpha_j, \lambda_q -\alpha_j)
A_{p,q}- (m_q +m_p) A_{p,q}$
\paragraph*{3.12.7} Let $\Omega$ be the non-commutative associative algebra 
generated by $\Omega (l,k), 1\leq l, k\leq n$.  Then for a fixed $k\neq l$.\\
$\{\Omega \cdot A_{k,l}\} =V^+_{\lambda-2\alpha j}$.

\paragraph*{(3.13)} We recall some well known facts from the representation 
theory of
$\mathfrak{sl}(2,\mathbb C)$ found in Humphreys' book [H].  Let
$\mathfrak{g}=\mathfrak{sl}(2,\mathbb C)$ with basis $x,y,h$  and $[x,y]=h$, 
$[h,x]=2x$ and $[h,y]=-2y$.    Let $m$ be a fixed positive integer and let 
$V(m)$ denote the finite dimensional irreducible highest weight module for 
$\mathfrak{sl}(2,\mathbb C)$ with highest weight vector $v$.
Then
\paragraph*{3.13.1} $hv=mv,\quad y^mv\neq 0,\quad y^{m+1}v=0.$
In Humphrey's book [H], Lemma 2.6.2 states

\paragraph*{(3.13.2)} $xy^a=y^ax+ay^{a-1}(h-a+1).$
for all $a\in \mathbb N$.

For $m,n\in\mathbb N$ one has the Clesbch-Gordan decomposition theorem

\paragraph*{(3.13.3)} $V(m)\otimes V(n)\cong V(m+n)\oplus V(m+n-2)\oplus 
\cdots \oplus  V(|m-n|)$
and this decomposition is multiplicity free.

Let us write down the highest weight vectors (up to a scalar) in this 
decomposition in terms of tensor products of weight vectors from $V(m)$ and 
$V(n)$.  Let $v_1$ and $v_2$ be the highest weight vectors of $V(m)$ and 
$V(n)$ respectively.   Then the highest weight vector of weight $w_1$ of 
weight $m+n-2l$ is a linear combination of the vectors $y^iv_1\otimes 
y^{l-i}v_2$ where $0\leq i\leq l$.
Then

\paragraph*{(3.13.4)} $w_l=\sum_{i=0}^la_iy^iv_1\otimes y^{l-i}v_2$
with $a_i\in\mathbb C$.  As $w_l$ is a highest weight vector we have
\begin{align*}
0&=xw_l=\sum_{i=0}^lxa_iy^iv_1\otimes y^{l-i}v_2+\sum_{i=0}^la_iy^lv_1\otimes
xy^{l-i}v_2 
\end{align*}
Thus one concludes
\paragraph*{(3.13.5)} $i(m-i+1)a_i+(l-i+1)(n-l+i)a_{i-1}=0$

for $1\leq i\leq l$. One can solve this recursion relation to obtain that the 
vector $(a_0,a_1,\cdots, a_l)$ is uniquely determined by just one of the 
coefficients say $a_0$ and each of the $a_i$ are nonzero.

Let $m>n>k$ and $m-n>k>0$.  Our goal now is to see how $V(m)\otimes 
V(n)\otimes V(k)$ decomposes and using the operators $\Omega(b_i,b_j)$ how to 
obtain a basis for all of the highest weight vectors in this tensor product.    
   The following matrix will explain this decomposition:
$$
\begin{pmatrix}
V(m+n+k) & V(m+n+k-2) &\cdots & V(m+n-k) \\
V(m+n+k-2) & V(m+n+k-4) &\cdots & V(m+n-k-2) \\
\vdots & \vdots &\ddots &\vdots  \\
V(m+n+k-2l) & V(m+n+k-2l-2) &\cdots & V(m+n-k-2l) \\
\vdots & \vdots &\ddots &\vdots  \\
V(m-n+k+2) & V(m-n+k) &\cdots & V(m-n-k+2) \\
V(m-n+k) & V(m-n+k-2) &\cdots & V(m-n-k) \\
\end{pmatrix}
$$
Let $V_{ij}:=V(m+n-2i+k-2j)$.  The matrix $(V_{ij})_{0\leq i\leq n,0\leq 
j\leq k}$ is the matrix above with $n+1$ rows and $k+1$ columns.  Notice that 
sum of the elements in the $(l+1)$-st row is nothing but the decomposition of
$V(m+n-2l)\otimes V(k)$.  Let $i+j=l$ with $0\leq i\leq n$, $0\leq j\leq k$, 
then $V_{ij}\cong V(m+n+k-2l)$.  

The sum of the first column is nothing but the decomposition of $V(m+k)
\otimes V(n)$. Similarily the sum of the last column is nothing but the
decomposition of $V(m-k) \otimes V(n)$.

The set of $V_{ij}$, $i+j=l$ is what we will call 
the {\it anti-diagonal} and they are all isomorphic.

Set $s_l=\min (l,k)$ and $d_l=s_l+1$  Then define

\paragraph*{(3.13.6)} $d_l'=\#\{(i,j)\,|\, i+j=l,0\leq i\leq n,0\le j\leq k\}$
It is easy to see that the following are true:
$$
d_l'=d_l,\quad 0\leq l\leq n,\qquad d_{n+i}'=k+1-i,\quad 0\leq i\leq k.
$$
Just for clarity we see $d_n'=k+1$ and $\min(n,k)=k$.

\paragraph*{(3.13.7)} Notice that the first entry of the $(l+1)$-row is the top 
component of $V(m+n-2l)\otimes V(k)$.  The highest weight vector  of this  
component is $w_l\otimes v_3$ where $v_3$ is the highest weight vector of 
$V(k)$.

Recall $w_l=\sum_ia_iy^iv_\otimes y^{l-i}v_2$ and each summand is nonzero.  
In particular $w_l\otimes v_3$ has $v_1\otimes y^lv_2\otimes v_3$ as a summand.

\paragraph*{(3.13.8)} Let $P_i =b_i$.  We claim $w_l\otimes v_3$ and 
$\Omega(b_2,b_3) (w_l\otimes v_3)$ are  linearly independent.
To prove the claim first note that $x(b_2)y(b_3)$ occurs in 
$\Omega(b_2,b_3)$.
Thus $\Omega(b_2,b_3)(w_l\otimes v_3)$ contains the term
$$
v_1\otimes y^{l-1}v_2\otimes yv_3
$$
and this term doesn't occur in $w_l\otimes v_3$.  Now the claim follows.

\paragraph*{(3.13.9)} The following set contains exactly $d_l$ linearly 
independent vectors for  $l\leq n$.
$$
\{ \Omega(b_2,b_3)^j(w_l\otimes v_3)\,|\, 0\leq j\leq s_l\}
$$
\paragraph*{Proof:}
 Note that by argument similar to the above we see that
 $$
 \Omega(b_2,b_3)^j(w_l\otimes v_3)
$$
contains the summand $v_1\otimes y^{l-j}v_2\otimes y^jv_3$ which doesn't
occur for lower $j$.
Thus the set consists of linearly independent vectors. This completes the
proof of the claim.

\paragraph*{(3.13.10)}Note that for $j>k$, the summands $v_1\otimes y^{l-j}
v_2\otimes y^jv_3$ is zero as $y^{k+1}v_3=0$. Similarly for $j>l$ that the 
term doesn't make sense.    Thus $j$ can go only up to $\min(l,k)=s_l$.

Up to now we have only worked with highest weight vectors 
$w_l\otimes v_3$, $0\leq l\leq n$. There are exactly $n+1$ highest weight 
vectors in the first column of the matrix.  By applying operators 
$\Omega(b_2,b_3)$ we get all of the highest weight vectors of the 
corresponding anti-diagonal.

\paragraph*{(3.13.11)}Now we will work with the highest weight vectors of the 
last row and prove that by applying $\Omega(b_2,b_3)$ repeatedly we can obtain
all 
other highest weight vectors.   Next consider the last row.  The first entry 
in the last row is taken care of.    We will only work with the second entry 
of the last row which is the representation
$$
V_{n,1}=V(m-n+k-2).
$$
Note that this module is the second component of $V(m-n)\otimes V(k)$.   
Since $w_n$ is the highest weight vector of $V(m-n)$ it is easy to see that
$$
z=(kyw_n\otimes v_3)-(m-n)w_n\otimes yv_3)
$$
is the highest weight vector of the second component of $V(m-n)\otimes V(k)$.   
Recalling the definition of $w_n$, we see that
$$
v_1\otimes y^nv_2\otimes yv_3
$$
is a nonzero summand of $z$ where we use the fact that $m>n$.

By applying $\Omega(b_2,b_3)^j$, with $j\leq k$ to $z$ we see that
$$
v_1\otimes y^{n-j}v_2\otimes y^jv_3
$$
is a summand of $\Omega(b_2,b_3)^jz$.   They are linearly independent and they  
are $k$ in number.  This is precisely the number of modules in the 
anti-diagonal as $d_{n+1}'=k$.
This argument breaks down for $j\geq k+1$ as $y^{k+1}v_3=0$.  Similarly the 
argument is valid for the other entries in the last row and we leave the 
details to the reader.

We will summarize the above results.  We have taken the highest weight vectors
of the first column and the last row.  Then we have applied our operators to 
 the highest weight vectors and obtained all other highest weight vectors.   

\paragraph*{4. Section}

\paragraph*{(4.1)} In this section we consider general linear algebra 
$\g = gl_N$ for a fixed positive  integer $N$.  Let $A$ be any commutative 
associated algebra with unit 1.  Then $\g\otimes A$ is a naturally Lie
algebra.  We will now define vectors in $U(\g\otimes A)$ which commutes
with $\g\cong \g \otimes 1$.  They are automatically central operators
on $\g\otimes A$ modules. Let
$$\{E_{ij}, \ 1 \leq i, \ j\leq N\}$$
be the standard basis with Lie bracket.
$$
[E_{ij}, \ E_{kl}]=\delta_{jk} E_{il} -\delta_{il}  E_{kj}
$$

\paragraph*{(4.2)} For a positive integer $k$ and $b_1, b_2,\ldots b_k$ in 
$A$, define 
$$
T_k(b_1,b_2,\ldots b_k) =\sum_{(i_1, i_2,\ldots i_k)} E_{i_1 i_2} (b_1)
 E_{i_2 i_3} (b_2)\ldots E_{i_ki_1} (b_k) 
$$
where $(i_1, i_2,\ldots, i_k)$ run over all possible indices.

Let $Z$ be the center of $U(\g)$. Then it is well known that $T_k(1,\ldots,1)
\in Z$ for all $k$.

\paragraph*{(4.3) Fact:}  It is a classical result of Harishchandra that 
$T_1(1,\ldots 1), \ldots, T_N(1,\ldots 1)$ generate $Z$ as an algebra.

\paragraph*{(4.4) Proposition:} Notation as above
$$
[T_k(b_1,\ldots, b_k), \g]=0
$$
\paragraph*{Proof} Clearly $T_k (b_1,\ldots, b_k) \in U (\g \otimes A)$. 
Note that $T_k(b_1, b_2,\ldots b_r + b^1_r, \ldots b_k) = T_k(b_1, b_2,\ldots,
b_r, \ldots, b_k) + T_k(b_1, b_2, \ldots b^1_r, \ldots b_k)$

Let $E_{j_1 j_2} \in \g$. For $ 1\leq r \leq k$.

Define
$$
B_r=\sum_{(i_l,\ldots i_k)} E_{i_1 i_2}(b_1)\ldots E_{i_{r-1}i_r}(b_{r-1}) 
E_{i_r j_2} (b_r) E_{j_1 i_{r+2}} (b_{r+1})\ldots E_{i_k i_1} (b_k)
$$
$$
B^1_r=\sum_{(i_l,\ldots i_k)} E_{i_1,i_2}(b_1) \ldots E_{i_{r-1}{j_2}} 
(b_{r-1}) 
E_{j_{1}i_{r+1}}(b_r)E_{i_{r+1}i_{r+2}} (b_{r+1}) \ldots E_{i_k i_1} (b_k)
$$

Now it is direct checking that 
$$
z=[T_k (b_1,\ldots, b_k), E_{j_1 j_2}]= \sum^k_{r=1} (B_r - B^1_r)
$$
Notice that $B_r=B^1_{r+1}$ for $1\leq r \leq k-1$ and $B_k =B^1_1$.

Now it is easy to see that $z=0$.  This proves the proposition.

\paragraph*{(4.5) Proposition (Remark 12 of [KS2]):} Let
$$
U(\g \otimes A)^\g =\{X \in U(\g \otimes A)\mid [\g, X]=0\}
$$
Then $T=U(\g\otimes A)^\g$.
\paragraph*{Proof} In [KS2] the proposition noted only for the polynomial algebra in one variable.
But the proof holds good for any commutative associative algebra $A$.

\paragraph*{(4.6)} In the rest of the section we take $A={\CC}
[t, t^{-1}]$.  Fix positive integer $n$. Let $\underline{\lambda}=
(\lambda_1,\ldots, \lambda_n)$ be dominant integral weights. Let 
$\underline{a}=(a_1,\ldots, a_n)$ be non-zero distinct complex numbers.  

For each $i$, let $V(\lambda_i)$ be an irreducible finite dimensional
 highest weight module with
highest weight vector $v_i$ for $\g$. Consider $V(\underline{\lambda}, 
\underline{a})=\displaystyle{\bigotimes^n_{i=1} V}(\lambda_i)$ is an 
irreducible evaluation
module for $\g\otimes A$. Recall from earlier section the polynomials
$P_1(t),\ldots, P_n(t)$ such that $\sum^n_{i=1} P_i(t) =1$

\paragraph*{(4.7) Remark.} $T_k(1\ldots, 1)$ which is a 
central operator and acts as scalar on every isotypic component of 
$V(\underline{\lambda}, \underline{a})$. But $T_k(1,1,\ldots, 1)$ splits into
 several operators.  $T_k(P_{i_1}(t),\ldots P_{i_k} (t))$ where each operator does
not act as scalars (most often). For clarity  we write one such operator.  Take
$n=4$ and consider $V(\lambda_1)\otimes \ldots \otimes V(\lambda_4)$.  Take
$k=3 $ and 
$$
T_3 (P_1(t), P_2(t), P_3(t)) w_1\otimes \ldots \otimes w_4=
\sum_{(i_1,i_2, i_3)}E_{i_1 i_2} w_1 \otimes E_{i_2i_3} w_2 \otimes E_{i_3 i_1}
 w_3 \otimes w_4
$$
Notice that there is no action on $w_4$. We believe such operators are
completely  new.

\paragraph*{(4.8) Theorem:} Let $V(\underline{\lambda}, \underline{a})= 
\oplus_\mu W(\mu)$ where $W(\mu)$ is an isotypic component.
Then each $W(\mu)$ is an irreducible modules for $T\otimes U(\g)$.
\paragraph*{Proof:} In this $V(\underline{\lambda}, \underline{a})$ is 
actually a moldule for $\oplus \g$ and $T \cong U(\oplus \g)^\g$ as operators.
Let 
$$
W(\mu)^+ =\{v\in W(\mu) \mid \g_+v=0\}
$$
Then  it is a well kown fact that $W(\mu)^+$ is an irreducible  module for
 $U(\oplus \g)^\g$ (double centralizer result, see [D, thm 9.1.12]).  
Now the theorem follows. 

\paragraph*{(4.9)} We will now extend the above results for the orthogonal
and symplectic Lie algebras. We will only sketch the results and leave the 
details to the reader.  We will follow closely Alexander Molev's book [Chapter
4 of [M]]. These Lie algebras can be treated as subalgebras of $gl_N$.  As
in the book we will number the rows and columns of $N\times N$ matrices by
the indices
$$
\{-k,\cdots, -1, 0,1,\cdots k\}
$$
if $N=2k+1$(orthogonal case) and by
$$
\{ -k,\cdots, -1,1,\cdots k\}
$$
if $N=2k$. (symplectic case).\\

Define 
$$
\theta_{ij}=\begin{cases}1 \ in \ the \ orthogonal \ case \\
(sign \ i\  sign \  j) \  in\  the \ symplectic  \ case
\end{cases}
$$
Note that $\theta^2_{ij} =1, \ \theta_{ij} =\theta_{ji}$ and 
$\theta_{ij} \ \theta_{jk} =\theta_{ik}$.

Let $\g_N$ denotes the one of these Lie-algebras.

Define
$$
F_{ij} =E_{ij} -\theta_{ij} E_{-j,-i}
$$
\paragraph*{(4.10)} The following is direct verification
$$
 [F_{ij}, F_{kl}]= \delta_{jk}  F_{il} -\delta_{il}  F_{kj}
+ \delta_{l,-j}  \theta_{ij} \ F_{k,-i}-\delta_{i,-k} \theta_{ij} F_{-j,l}
$$
As earlier let $A$ be a commutative associative algebra with unit. We will now 
construct central operators for the Lie algebra $\g_N \otimes A$. Let $b_1,
b_2,\ldots b_k\in A$ and $k > 0$.

\paragraph*{(4.11) Define}
$$
S_k (b_1,\ldots, b_k)= \sum_{(i_1,\ldots i_k)} F_{i_1 i_2} (b_1)\ldots
F_{i_k i_1}(b_k)
$$
It is direct checking that the above operators are central.  

\paragraph*{(4.12) Remark:} Results similar to $gl_N$ also holds for type $B$ and $C$,
\paragraph*{5. Section Spanning set for $T$.}

In this section we take $\g =gl_N$ for some positive integer $N$.  Let $A$
be commutative associative algebra with unit.  Recall that $T$ is a non-commutative associative subalgebra of $U(\g \otimes A)$ generated by $T_k(a_1,
a_2,\ldots, a_k)$.  See 4.2. The purpose of this section is to give a spanning
set for $T$ and avoiding products.

\paragraph*{(5.1)} Throughout this section we will be dealing with following
finite sets
$$
S=\{ (i_j,i_k) \mid i_j\ and \ i_k \ are \ variables, \ j \ 
and \ k \ denote \ some \ positive \ integers\} 
$$

Let $(i_j, i_k), (i_m,i_n) \in S$. Then they are said to be connected
if $i_k=i_m$ or $i_j=i_n\cdot i_j$ is called the start point and $i_k$ is 
called end point of $(i_j, i_k)$.

\paragraph*{5.2. Definition.}

A finite set $S$ is called circuit of\\
(1) Each start point and end point occurs only once\\
(2) If $i_k$ is start point (resp. end point)
then it is also occurs as  end point (resp. start point).\\
(3) $S$ is connected in the sense if $(i_j, i_k), (i_m, i_n)\in S$
then there eixsts a sequence of elements in $S$ starting with $(i_j,i_k)$ and ending with $(i_m,i_n)$ and consecutive elements  are connected.

\paragraph*{Remark.}
The first two conditions imply that $S$ is union of circuits.

\paragraph*{(5.3) Example.}
Let $n$ be a positive integer.  Then
$$
S=\{(i_1,i_2), (i_2,i_3) \ldots (i_n,i_1)\}
$$
is a circuit.  In fact any circuit is of this form after rearranging the 
indices. We denote $\mid\!\!S\!\!\mid$ the number of elements of $S$.

\paragraph*{(5.4) Definition.}
 Let $S=\cup S_j$ be disjoint union of circuits.  Suppose $(i_{k-1},i_k), (i_l,
i_{l+1})\in S$.  Then define $S_{k,l}$ be the set of elements of $S$ excluding
$(i_{k-1}, i_k), (i_l,i_{l+1})$ and including $(i_{k-1},i_{l+1})$.  Further
replace $i_l$ with $i_k$ everywhere. Clearly $\mid\!\!S\!\!\mid =\mid\! 
S_{k,l}\!\mid +1$.

\paragraph*{(5.5) Lemma.}

$S_{k,l}$ is union of circuits

\paragraph*{Proof:}
We will first assume both

\paragraph*{(5.5.1)} $(i_{k-1},i_k),(i_l,i_{l+1})\in S_j$ for some $j$.\\
We further assume $S_j=\{(i_1,i_2),(i_2,i_3),\ldots (i_{m_j},i_1)\}$ and
$\mid\!\! S_j\!\!\mid =m_j$.  It is easy to check the Lemma when $m_j$ is very small. 
Thus can assume $m_j \geq 4$.

Suppose $k=l$ then the Lemma is obivious.  Assume $k < l$.  Then clearly
$$
S_{k,l} =\{ (i_1,i_2)\ldots (i_{k-2}, i_{k-1}), (i_{k-1}, i_{l+1}),
(i_{l+1}, i_{l+2}),\ldots (i_{m_j}, i_1)\}\cup\\
$$
$$ \{(i_k, i_{k+1}),\ldots (i_{l-1}, i_k)\}
$$
Note that $i_l$ is replaced by $i_k$. Thus $S_{k,l}$ is union of circuits.  Now assume $l<k$.  But one can order $S_j$ in such a way that $l$ occurs 
after $k$.  Then the Lemma follows from earlier case.

Now suppose element in (5.5.1) occurs in different circuits.  We can assume
 $S=S_1 \cup S_2$ and 
$$
\begin{array}{lll}
S_1 &=& \{ (i_1,i_2),\ldots (i_{m_1},i_1)\}\\
S_2 &=& \{ (j_1, j_2),\ldots (j_{m_2}, j_1)\}\\
\end{array}
$$
and $(i_{m_1},i_1), (j_1, j_2)$ are elements in (5.5.1).

Now the Lemma is obivious noting 
\paragraph*{(5.5.2)} $S_{k,l} =\{(i_1,i_2),\ldots(i_{m_1}, j_2),(j_2,j_3),
\ldots (j_{m_2}, i_1)\}$ which is a single circuit.

\paragraph*{(5.6)} We wil now defne certain twisted product and prove that
they are central.

Let $S=\cup S_j$ be disjoint union of circuits and let $\mid\!\!S_j\!\!\mid 
=m_j$ and  $m=\sum m_j =\mid\!\! S\!\!\mid$. Let

\paragraph{(5.6.1)} $\underline{a}_m =(a_1,a_2,\ldots a_m) \in \oplus A=A_m 
(m$ copies)

Let $(i_{j_1},i_{k_1}),\ldots (i_{j_m}, i_{k_m})$ be some order of element in 
$S$. We will denote this permutation of $S$
by $\sigma$.  Every  circuit has a natural  order (not unique) in the sense 
that the consecutive elements are connected.  When the order of $S_j$ in 
natural  we denote the permutation by $I_d$.  Note that if $(i_{j_t}, i_{k_t})
\in S_j$ then $j_t +1 =k_t$ (read  mod $m_j$).

Define

\paragraph{(5.6.2)} $T_m (S,\sigma, \underline{a}_m)
=\sum_{(j_1,j_2,\ldots j_m)}
 E_{i_{j_1} i_{k_1}} (a_1)\ldots E_{i_{j_m} i_{k_m}}(a_m) 
$\\

where the summation runs over all possible indices from 1 to N. We will say the
order of the above operator is $m$.

\paragraph{(5.6.3)} (a) Operator of the above form are called twisted 
product.\\
(b) Recall that $T_k(a_1,a_2,\ldots a_k)=T_k(S^{'}, id, \underline{a}_k)\in T$ 
where \\
$S'=\{(i_1,i_2),\ldots (i_k,i_1)\}$. Product of such operators are 
called straight product.  As mentioned earlier natural order is not unique 
but the corresponding operator is same.

\paragraph{(5.7) Proposition:} The operator $T_m(S,\sigma, \underline{a}_m)$ is
central.  We need the following

\paragraph*{(5.8) Lemma}
$T_m (s,\sigma, \underline{a}_m) =\prod T_{m_i}(S_i, Id,
\underline{b}_{m_i})$ + lower order twisted operators for some 
$\underline{b}_{m_i} \in A_{m_i}$.

\paragraph{Proof} The proof is very simple. By interchanding consecutive $E's$
in the product of $T_m(S,\sigma, \underline{a}_m)$ we can get to the first 
term of the right hand side.  Every time we interchange two $E's$ we get two 
additional twisted product but of lower order. We will explain this in more 
detail.  Let us say we interchange $(i_{k-1}, i_k)$ $(i_j, i_{j+1})$ which are 
consecutive  entries.
$$
\begin{array}{lll}
&T_m(S,\sigma, \underline{a}_m)=\sum E_{i_{j_1} i_{k_1}} (a_1)\ldots 
E_{i_{k-1} i_k} (a^{'}) E_{i_j i_{j+1}} (a^{''})\ldots E_{i_{j_m} i_{k_m}}
 (a_m)\\
&=\sum E_{i_{j_1} i_{k_1}} (a_1)\ldots E_{i_j i_{j+1}}(a^{''}) E_{i_{k-1} i_k}
(a^{'}) \ldots E_{i_{j_m} i_{k_m}} (a_m)\\
&+\sum_{i_k=i_j} E_{i_{j_1} i_{k_1}} (a_1) \ldots E_{i_{k-1} i_{j+1}} 
(a^{'} a^{''}) \ldots E_{i_{j_m}i_{k_m}} (a_m)\\
&-\sum_{i_{j+1}=i_{k-1}} E_{i_{j_1}} i_{k_1} (a_1)\ldots E_{i_j i_k} (a^{'} 
a^{''}) \ldots E_{i_{j_m} i_{k_m}} (a_m)
\end{array}
$$
Notice the sets corresponding to the three operators on the right hand side 
are $S, S_{k,l}, S_{l+1,k-1}$.  They are all union of circuits by Lemma 5.5. 
Further $\mid\!\! S_{k,l}\!\!\mid =\mid\!\! S_{l+1,k-1}\!\!\mid =\mid\!\! 
S\!\!\mid-1$  By repeating 
this process several times we complete the proof of the Lemma.

\paragraph{(5.9) Corollary:} $T_m (S,\sigma, \underline{a}_m)$ equals to sum 
of straight products. (See 5.6.3(b) for definition).  Just apply above Lemma 
for lower order operators.

\paragraph{Proof of Proposition 5.7} Since straight products are central the
proposition follows.

Let $\stackrel{\sim}{T}$ be linear span of $T_m (S,\sigma, \underline{a}_m), 
m\in \NN, S$ is any single circuit such that $\mid\!\! S\!\!\mid =m, 
\sigma$ is any 
order of $S$ and for all  $\underline{a}_m \in A_m$

\paragraph{(5.10) Theorem:} $\widetilde{T} =T$.

\paragraph{Proof} By definition $T$ contains all straight products.  By 
corollary 5.9 each $T_m(S,\sigma, \underline{a}_m)$ is sum of straight 
products and hence $\widetilde{T} \subseteq T$.

\paragraph*{Claim(1)} $T_m (S,\sigma, \underline{a}_m) T_l (b_1,b_2,\ldots 
b_l)\in \widetilde{T}$\\

(2) $T_l (b_1,b_2,\ldots b_l) T_m (S,\sigma, \underline{a}_m) \in 
\widetilde{T} . \ l \in \NN, b_i \in A, \ \underline{a}_m \in A_m$.

We will first complete the proof of the Theorem by assuming the claim.  By
claim it follows that the straight product of two operators is in 
$\widetilde{T}$. Again by claim we see straight product of three operator is 
in $\widetilde{T}$. Similarly any straight product is in $\widetilde{T}$. 

This proves $T\subseteq \widetilde{T}$.  This completes the proof of Theorem.

\paragraph*{Proof of the Claim 1} $S$ comes with some order and let that order 
be $(i_{j_1}, i_{j_1+1}),\ldots (i_{j_m}, i_{j_{m}+1})$. For 
$l\geq 0$ consider $\widetilde{S}$ with the order $\widetilde{S}=$
$$\{ (i_{j_1}, i_{j_1+1}),\ldots(i_{j_{m-1}}, i_{j_{m-1}+1}), 
(i_{j_m}, i_{j_m+2}), (i_{j_{m+2}}, j_1), (j_{l+1}, i_{j_m+1}), (j_1,j_2),
\ldots (j_l, j_{l+1})\}.
$$
We have obtained $\widetilde{S}$ from $S$ by deleting $(i_{j_m}, i_{j_m+1})$
and adding \\
$(i_{j_m}, i_{j+2}), (i_{j_{m+2}}, j_1), (j_{l+1},
 i_{j_m+1}), (j_1,j_2),\ldots (j_l, j_{l+1})$.

It is easy to check that $\widetilde{S}$ is a single circuit.  It comes with 
an order and  denote it by $\sigma^1$. From the definition it follows that
$T_{m^1}(\widetilde{S}, \sigma^1, \underline{d}_{m^1})\in \widetilde{T}$
where $\underline{d}_{m^1}=(\underline{a}_m, 1, b_{l+1}, b_1,\ldots b_l)$. Note
that $\mid\!\!\widetilde{S}\!\!\mid=m^1 =m+l+2$. Notice that 
$E_{i_{j_{m+2}} j_1}(a_m) E_{j_{l+1} i_{j_m+1}}(1)$ occurs in 
$T_m(\widetilde{S}, \sigma^1, \underline{d}_{m^1})$.  

As in the earlier argument we interchange these two term and we obtain the 
following equation.

\paragraph*{(5.10.1)} $T_{m^1}(\widetilde{S},\sigma^1, \underline{d}_{m^1})
-T_{m^1} (\widetilde{S}, \sigma^{''}, \underline{d}_{m^{''}})=N T_{m+1} (S_1,
\sigma_1, \underline{d}_{m+1}) T_l (S_2, 1d, \underline{b}_l) -N T_m
(S,\sigma, \underline{a}_m). T_{l+1} (S_3, 1d, \underline{b^{'}}_{l+1})$
where $\underline{d}_{m^{''}}$ is obtained from $\underline{d}_{m'}$ by interchanding 1 and $b_{l+1}$. 

$$
\begin{array}{lll}
\underline{d}_{m+1} &=& (a_1, a_2, \ldots, a_m, b_{l+1})\\
\underline{b}_l &=& (b_1,b_2,\ldots, b_l)\\
\underline{b}^{'}_{l+1} &=& (b_{l+1}, b_1, b_2,\ldots, b_l)\\[4mm]
S_1&=&\{(i_{j_1},i_{j_1+1}),\ldots (i_{j_{m-1}}, i_{j_{m-1}+1}),(i_{j_m}, 
i_{j_{m+2}}),(i_{j_{m+2}} i_{j_m+1})\}\\
S_2 &=& \{(j_1,j_2),\ldots, (j_l, j_1)\}\\
S_3& =&\{(j_{l+1}, j_1), (j_1,j_2)\ldots (j_l, j_{l+1})\}
\end{array}
$$
The order $\sigma^{''}$ is obtained from $\sigma^{'}$ by interchanging 
$(i_{j_{m+2}}, j_1)$ and $(j_{l+1}, i_{j_m+1})$.  The order $\sigma_1$ is 
given in definition of $S_1$.  The order $\sigma$ is the one we started with.
Since the terms in left hand side are in $\widetilde{T}$, the difference of 
terms in the right hand side is also in  $\widetilde{T}$.
Suppose $l=0$ then \\
$\widetilde{S} =\{(i_{j_1}, i_{j_1+1}),\ldots
(i_{j_{m-1}},i_{j_{m-1}+1}), 
 (i_{j_m}, i_{j_{m+2}}), (i_{j_{m+2}},j_1), (j_1, i_{j_m+1})$

The RHS is $N T_{m+1} (\widetilde{S}, \sigma^{'}, \underline{d}_{m+1}) Id -N T_m
(\widetilde{S}, \sigma, \underline{a}_m) T_1 (S_3, Id, b^{'}_1)$

>From this we conclude that $\widetilde{T}$ is closed  under right 
multiplication by $T_1(S_3, Id, b^{'}_1) =T_1 (b^{'})$. Now using induction on
 $l$ and by (5.10.1) we see that $\widetilde{T}$ is closed under multiplication by
$T_l (b_1,\ldots, b_l)$.  This proves Claim 1.  Claim 2 is similar. 
This completes the Proof of the Theorem.

\paragraph*{(5.11 Remark)} Let $T_m (S^1, \sigma^1, \underline{a}_m), T_n (S^2,
\sigma^2, \underline{b}_n) \in \widetilde{T} =T$. We are assuming both $S^1$ 
and $S^2$ are single circuits.
$$
[T_m(S^1, \sigma^1, \underline{a}_m), T_n (S^2, \sigma^2, \underline{b}_n)]
$$
$$
=\sum_{\substack{1\leq k\leq m\\ 1\leq l \leq n}}(T_{m+n-1}(S_{k,l}, 
\sigma_{k,l},  \underline{d}^{k,l}_{m+n-1})-T_{m+n-1} (S_{l+1,k-1}, 
\sigma^1_{k,l}, \underline{d}^{k,l}_{m+n-1})
$$

where $S_{k,l}$ is a single circuit.  $\sigma_{k,l}, \sigma^1_{k,l}$ are some 
pernutations.  $\underline{d}^{k,l}_{m+n-1} \in A_{m+n-1}. (l+1, k-1$ read
mod $n$ and $m$).

It is very elmentary to see the  remark.  We will explain $S_{k,l}$
Let
$$
\begin{array}{lll}
S&=&S^1\cup S^2 \ where\\
S^1 &=&\{ (i_1,i_2),\ldots, (i_m, i_1)\}\\
S^2&=& \{ (j_1, j_2), \ldots, (j_n, j_1)\}
\end{array}
$$
Let $(i_{k-1}, i_k) \in S^1, (j_l, j_{l+1})\in S^2$.
Then $S_{k,l}$ is defined in 5.4. Similarly one can define $S_{l+1, k-1}$.
In Lemma 5.5. we noted that each $S_{k,l}$ is a single circuit.

\newpage
\paragraph*{References}
\begin{enumerate}

\item[{[D]}] J. Dixmier, Enveloping algebras, North Holland, 1977.

\item[{[EB]}] S. Eswara Rao and Punita Batra, Classification of irreducible
integrable highest weight modules  for current Kac-Moody Lie algebra,
Journal of Algebra and its Applications, Vol 16, No.5 (2017), 1750123.

\item[{[EF]}] S. Eswara Rao and V. Futorny, Representations of loop Kac-Moody
Lie algebras, Comm. Algebra  41(2013), No. 10, 3775-3792.


\item[{[EFS]}] S. Eswara Rao, V. Futorny and  Sachin Sharma, Weyl modules 
associated to Kac-Moody lie algebras, Communications of Algebra 44(2016), No. 12, 5045-5057.

\item[{[E1]}] S. Eswara Rao, On representations of loop algebras, Comm.
Algebras, 21(1993), No.6, 2131-2153. 

\item[{[E2]}] S. Eswara Rao, Classification of irreducible integrable modules 
for multi-loop algebras with finite  dimensional weight spaces, Journal of
Algebra, 246 (2001), 215-225.

\item[{[E3]}] S. Eswara Rao, Classification of irreducible integrable modules
for toroidal Lie algebras with finite dimensional weight spaces, Journal of 
Algebra, 277 (2004), 318-348.

\item[{[G]}]  Werner Greub, Linear algebra, Springer-Verlag, New York-Berlin
fourth edition, 1975, Graduate Texts in Mathematics, No.23.
 
\item[{[GE]}] I.M. Gelfand, Center of the infinitesimal group ring,. Math.
Sboronik 26 (1950), 103-112 (Russian) English tranl. in I.M. Gelfand, 
`Collected papers', Vol. II: Springer-Verlag, 1988, PP 22-30.

\item[{[H]}] J.E. Humphreys, Introduction to Lie algebras and representation
theory, Vol. 9 of Graduate Texts in Mathematics, Springer Verlag, New York 1978.

\item[{[K]}] V.G. Kac, Infinite dimensional Lie algebras, Cambridge university
press, Cambridge, 3rd edition 1990.

\item[{[KS1]}] S. Kumar, Tensor product decomposition, proceedings of 
International Congress of Mathematicians, Hyderabad, India Vol. III, 1226-1261,
New Delhi, 2010.

\item[{[KS2]}] S. Kumar, A complete set of intervwiners for arbitary tensor 
product representation via current algebra, arxiv: 1607.06115 (to be 
published in Transformation Groups).

\item[{[M]}] A. Molev, Yangians and classical Lie algebra, SURU-143, AMS(2007).

\item[{[MU]}] I.M. Musson, Lie superalgebras and enveloping algebras, Vol.131,
 Graduate Studies in Mathematics, AMS.

\item[{[NS]}] E. Neher and A. Savage, A survey of equivariant map algebras 
with open problems, Recent developments in algebraic and combinatorial 
aspects of representation theory, Contemp. Math. 602, 165-182, AMS, Providence
RI, 2013.

\item[{[PRV]}] K.R. Parthasarathy, R. Ranga Rao and V.S. Varadarajan,
Representation of complex semi-simple Lie groups and Lie algebras, Ann. of
Math. 85 (1967), 383-429

\item[{[S]}] A. Savage, Equivariant map superalgebras, Math. $Z$, 277 
(2014) No.1-2, 373-399.

\item[{[SN]}] N. Sthanumoorthy, Introduction to finite and infinite 
dimensional lie (super) algebras, Academic Press (2016).

\item[{[W]}] M. Wakimoto, Infinite dimensiolnal Lie algebras, Translation of 
Mathematical Monographs, AMS, Vol 195, 2001.
\end{enumerate}

\noindent
School of Mathematics\\
Tata Institute of Fundamental Research\\
email: senapati@math.tifr.res.in, sena98672@gmail.com\\[5mm]

\noindent
Added in the Proof:  Although we have incorporated Shrawan Kumar's
[KS2] results, his preprint appeared much much later than our preprint
in arxiv.org.

\noindent
{\bf Acknowledgements:}\\
\noindent
I would like to thank Maria Gorelik for some helpful discussions.

\end{document}